\documentclass[11pt]{amsart}
\usepackage{amsfonts}
\usepackage{latexsym}
\usepackage{amsmath}
\usepackage{amssymb}
\newcommand{\be}{\begin{equation}} \newcommand{\ee}{\end{equation}}
\newcommand{\bea}{\begin{eqnarray}} \newcommand{\eea}{\end{eqnarray}}
\newcommand{\bean}{\begin{eqnarray*}}
  \newcommand{\eean}{\end{eqnarray*}}
\newcommand{\brray}{\begin{array}} \newcommand{\erray}{\end{array}}
\newcommand{\ben}{\begin{equation}{nonumber}}
  \newcommand{\een}{\end{equation}{nonumber}}
\newcommand{\newsection}[1]{\setcounter{equation}{0}
  \setcounter{dfn}{0}
\section{#1}}

\newtheorem{dfn}{Definition}[section]
\newtheorem{thm}[dfn]{Theorem}
\newtheorem{lmma}[dfn]{Lemma}
\newtheorem{ppsn}[dfn]{Proposition}
\newtheorem{crlre}[dfn]{Corollary}
\newtheorem{xmpl}[dfn]{Example}
\newtheorem{rmrk}[dfn]{Remark}

\newcommand{\bdfn}{\begin{dfn}} \newcommand{\bthm}{\begin{thm}}
    \newcommand{\blmma}{\begin{lmma}}
      \newcommand{\bppsn}{\begin{ppsn}}
        \newcommand{\bcrlre}{\begin{crlre}}
          \newcommand{\bxmpl}{\begin{xmpl}}
            \newcommand{\brmrk}{\begin{rmrk}}

              \newcommand{\edfn}{\end{dfn}}
            \newcommand{\ethm}{\end{thm}}
          \newcommand{\elmma}{\end{lmma}}
        \newcommand{\eppsn}{\end{ppsn}}
      \newcommand{\ecrlre}{\end{crlre}}
    \newcommand{\exmpl}{\end{xmpl}} \newcommand{\ermrk}{\end{rmrk}}
\newcommand{\ots}{\otimes}

\newcommand{\IC}{{\mathbb C}}

 \newcommand{\IN}{{\mathbb N}}
 
 \newcommand{\IR}{{\mathbb R}}
 \newcommand{\IT}{{\mathbb T}}

 \newcommand{\IZ}{{\mathbb Z}}

 \newcommand{\del}{\partial}

\newcommand{\cla}{{\mathcal A}} \newcommand{\clb}{{\mathcal B}}

 \newcommand{\clh}{{\mathcal H}}
 
\newcommand{\clk}{{\mathcal K}}

  \def 
\bbE {\mbox{\boldmath $E$}}                     

  \def 
\bbe {\mbox{\boldmath $e$}}

 \def\a*{{\mathcal A}_{h,*}} \def\B{{\mathcal B}(h)}
\def\B1{{\mathcal B}_1(h)} \def\b{{\mathcal B}^{s. a. }(h)} \def\b1{{\mathcal
    B}^{s. a. }_1(h)}

\newcommand{\raro}{\rightarrow}

\def \qed {  \mbox{}\hfill $\Box$\vspace{1ex}}

\begin{document}
\title[From C*algebra extensions to CQMS]{From C*algebra extensions to CQMS, ${SU}_q(2)$, Podles
sphere and other examples}

\author[Partha Sarathi Chakraborty]{Partha Sarathi Chakraborty}
\address[Partha Sarathi Chakraborty]{Stat--Math Division,
Indian Statistical Institute, 203, B. T. Road, Kolkata 700 108,
India.}
\email{\tt parthasc\_r@isical.ac.in,
partha\_sarathi\_c@hotmail.com}
\date{}
\begin{abstract}
We construct compact quantum metric spaces (CQMS) starting with
some C*algebra extension with a positive splitting. As special
cases we discuss the case of Toeplitz algebra, quantum SU(2) and
Podles sphere.
\end{abstract}
 \maketitle
\pagebreak
\section {Introduction}
In noncommutative geometry, the natural way to specify a metric is
by a ``Lipschitz seminorm''. This idea was first suggested by
Connes~(\cite{C5}), and developed further in \cite{C7}. Connes
pointed out~(\cite{C5}, \cite{C7}) that from a Lipschitz seminorm
one obtains in a simple way an ordinary metric on the state space
of a $C^*$-algebra. A natural question in this context is when
does this metric topology coincides with the weak* topology. In
his search for an answer to this question, Rieffel~(\cite{RI7},
\cite{RI8}, \cite{RI9}) has identified a larger class of spaces,
namely order unit spaces on which one can answer these questions.
He has introduced the concept of Compact Quantum Metric Spaces
(CQMS) as a generalization of compact metric spaces, and
used~(\cite{RI9}) this new concept for rigorous study of
convergence questions of algebras much in the spirit of
Gromov-Hausdorff convergence. One natural question in this regard
is are there plenty of CQMS floating around? Rieffel (\cite{RI7},
\cite{RI8}) has given some general principles for the construction
of CQMS. In \cite{CH} we exploited one of his principles to
construct CQMS. In fact Rieffel   has shown  (\cite{RI10}) that
indeed there are many examples. But in concrete C*-algebras one
would like to have more explicit description of these structures.
Our objective here is construction of CQMS out of quantum SU(2)
and Podles spheres. To achieve that we make a slightly general
construction and produce CQMS starting from C*-algebra extensions.
Organization of the paper is as follows. In the next section we
recall the basics of CQMS. In section 3 the basic construction has
been described. In the final section we employ the principle
developed in section 3 to special cases.
\section { Compact Quantum Metric Space: Preliminaries}
We recall some of the definitions from \cite{RI9}.
 \bdfn  \rm
An order unit space is a real partially ordered vector space $A$  with a distinguished element $e$, the order unit  satisfying \\
(i) ( Order Unit property ) For each $ a \in A$ there is an $ r \in \IR$ such that $ a \le r e $. \\
(ii) ( The Archimidean property ) If $ a \in A$ and if $ a \le r e $ for all $ r \in \IR$ with $ r \ge 0$, then $ a \le 0$.\\
\edfn \brmrk  \rm
 The following prescription defines a norm on an order
unit space.
\[ \| a \| = \inf \{ r \in \IR | -re \le a \le re \} \]
\ermrk \bdfn \rm  By a state of an order unit space $(A,e)$ we
mean a $\mu \in A^\prime$, the dual of $(A, \| \cdot \|)$ such
that $ \mu ( e) =1 = {\| \mu \|}^\prime$. Here $ { \| \cdot
\|}^\prime$ stands for the dual norm on  $A^\prime$. Collection of
states on $(A,e)$ is denoted by $S(A)$. \edfn \brmrk   \rm
 States are
automatically positive. \ermrk \bxmpl   \rm Motivating example of
the above concept is the real subspace of selfadjoint elements in
a C*-algebra with the  order structure inherited from the
C*-algebra. \exmpl \bdfn   \rm
Let $(A,e)$ be an order unit space. By a Lip norm on $A$ we mean a seminorm $L,$ on $A$ such that \\
(i) For $ a \in A$, we have $L(a)=0$ iff $ a \in \IR e$; \\
(ii) The topology on $S(A)$ coming from the metric  \[\rho_L (\mu,
\nu ) =sup \{ | \mu (a) -  \nu (a) | : L(a) \le 1 \} \] is the
$w^*$ topology. \edfn \bdfn   \rm A compact quantum metric space
is a pair $(A,L)$ consisting of an order unit space $A$ and a Lip
norm $L$ defined on it. \edfn

The following theorem of Rieffel will be of crucial importance.
\bthm[{\bf Theorem 4.5 of \cite{RI9}}] \rm  \label{P:thm:7.4}
 Let $L$ be a seminorm on the order unit space $A$ such that $L(a)=0 $ iff $a \in \IR e$. Then $\rho_L$ gives $S(A)$ the $w^*$-topology exactly
 if\\
(i) $(A,L)$ has finite radius, i.e,   $ \rho_L (\mu, \nu ) \le C
\mbox   { for all }  \mu, \nu   \in S(A)$ for some constant $C$,
and \\
(ii) $\clb_1= \{ a | L(a) \le 1,  \mbox{ and } \| a \| \le 1 \}$
is totally bounded in $A$ for $ \| \cdot \| $. \ethm
\newsection{Extensions to CQMS}

In this section we describe the general principle of construction
of CQMS from certain $C^*$-algebra extensions. Let $\cla$ be a
unital $C^*$-algebra. Fix a faithful representation $\cla
\subseteq \clb(\clh)$. Suppose we have a dense order unit space
$Lip(\cla) \subseteq \cla_{s.a}$, where $\cla_{s.a}$ denotes the
real partially ordered subset of selfadjoint elements in $\cla$.
Let $L$ be a Lip norm on $Lip(\cla)$ such that $((Lip(\cla),I),L)$
is a CQMS. Let $\nu$ be a state on $\cla$, then define $
\widetilde {\cla_\nu} $ to be the collection of $ ((a_{ij})) \in
\clk ( l^2(\IN)) \otimes \cla $ such that  (i)  $a_{ij}\in
Lip(\cla) $, (ii) $ a_{ij}=a_{ji} $, (iii) $sup_{ i \ge 1, j \ge
1} {(i+j)}^k ( L(a_{ij})+ |\nu(a_{ij}) | )< \infty \; \forall k $.
Clearly $\cla_\nu:=\widetilde {\cla_\nu}\oplus \IR I$, where $I$
is the identity on $\clb(\l^2(\IN) \otimes \clh)$ is an order unit
space. Define $L_k : \cla_\nu \raro \IR_{+} $  by  $ L_k (I) =0,$
\[
L_k((a_{ij}))= sup_{ i \ge 1, j \ge 1} {(i+j)}^k ( L(a_{ij})+
|\nu(a_{ij}) | ).
\]
\blmma Let $ d= $ diameter of $((Lip(\cla),I),L)$. Then for   a
``Lipschitz function" $a \in Lip(\cla)$ one has  $ \| a \| \leq (
L(a) + | \nu(a) |)(1+d). $ \elmma {\it  Proof:} Let $\mu$ be an
arbitrary state on $\cla$. Then using $\sup \{ |\mu(a)-\nu(a)| :
L(a) \le 1 \} \le d$ we get, \bean | \mu(a) |  & \leq  & | \mu(a)
- \nu(a) | + | \nu(a) | \cr
  & \leq & L(a )d  + | \nu(a) | \cr
& \leq & ( L(a) + | \nu(a) |)(1+d). \eean \qed \blmma
\label{M:lem:2.2} There exists a constant $C>0$ such that for $
((a_{ij})) \in  \widetilde {\cla_\nu},$
\[
\| ((a_{ij})) \| \leq C L_2((a_{ij})).
\]
\elmma {\it Proof:} Let ${\{ e_i \}}_{i \geq 1} $ be the canonical
orthonormal basis for $ l^2( \IN)$. Let $\sum \lambda_i e_i
\otimes u_i $ and  $\sum \mu_i e_i \otimes v_i $ be two generic
elements in $ l^2( \IN) \otimes \clh $. Here $ u_i, v_i \in \clh $
are unit vectors.  Then clearly ${ \| \sum \lambda_i e_i \otimes
u_i \| }^2 = {\sum | \lambda_i  | }^2,   { \| \sum \mu_i e_i
\otimes u_i \| }^2 = {\sum | \mu_i  | }^2. $ Now observe that
\bean | \langle \sum \lambda_i e_i \otimes u_i, ((a_{ij}))  \sum
\mu_j e_j \otimes v_j \rangle | & \leq & \sum | \lambda_i| | \mu_j
| | \langle u_i, a_{ij} v_j \rangle | \cr & \leq  & \sum |
\lambda_i| | \mu_j |   ( L(a_{ij}) + | \nu(a_{ij}) |)(1+d) \cr
  & \leq & (1+d)  \sum | \lambda_i| | \mu_j |  \frac {L_2 ((a_{ij})) }{ i j} \cr
& \leq & L_2 ((a_{ij})) (1+d ) \sum_{n=1}^\infty \frac {1} { n^2}
\sqrt{ \sum {|\lambda_i|}^2} \sqrt{ \sum {|\mu_i|}^2}. \eean This
proves the lemma with $C=(1+d ) \sum_{n=1}^\infty \frac {1} {
n^2}$. \qed \blmma Let $ \clb_1 = \{ a \in \cla_\nu | L_k (a )
\leq 1, \| a \| \leq 1 \}.$  Then  $ \clb_1$ is totally bounded in
norm for $k > 2 $. \elmma {\it Proof:} Let $ \epsilon > 0 $ be
given. Choose $N$ such that $ {( \frac {1}{N})}^{k-2} < \epsilon$.
For $G =(( g_{ij})) \in \cla_\nu$ let $ P_N(G) \in  \clk (
l^2(\IN)) \otimes \cla  $ be the element given by \\
$${P_N(G)}_{ij}= \begin{cases}  g_{ij} & \text{ for $ i,j \leq N,$}  \cr 0 &
\text{ otherwise.} \cr \end{cases}$$ Now observe that \bean L_k( G
- P_N(G)) & = & sup_{ i \geq N {\mbox or } j \geq N } {(i+j)}^k (
L (g_{ij}) + | \nu(g_{ij}) | ) \cr & \geq & N^{k-2} sup_{ i \geq N
{\mbox or } j \geq N } {(i+j)}^2  ( L (g_{ij}) + | \nu(g_{ij}) | )
\cr & = & N^{k-2}  L_2( G - P_N(G)). \eean Note that for $G \in
\clb_1 , L_k( G - P_N(G)) \leq 1$, therefore \bean \| G- P_N(G) \|
& \leq & C L_2( G - P_N(G)) \cr & \leq & C {N}^{-(k-2)}  L_k( G -
P_N(G)) < C \epsilon. \eean Here the constant $C$ is the one
obtained in the previous lemma. Note $C$ does not depend on $N$.
By theorem~\ref{P:thm:7.4} there exists $ N \times N$ matrices $
(( a_{ij}^{(r)} )) \in M_N (\cla),$   for   $r=1, \ldots, l$ such
that for any  $ N \times N$ matrix $ (( a_{ij} )) \in \clb_1 ,$
there exists $r$ satisfying $\|  (( a_{ij} )) -(( a_{ij}^{(r)} ))
\| < \epsilon $. Now for $G \in \clb_1 $, get $(( a_{ij}^{(r)} ))
$ such that $\|  P_N(G) -(( a_{ij}^{(r)} )) \| < \epsilon $. Then,
$$
\| G -(( a_{ij}^{(r)} )) \|   \leq   \|  G - P_N(G)\| + \epsilon
\leq (1+ C) \epsilon. $$ This completes the proof. \qed \bthm
$((\cla_\nu,I),L_k)$ is a compact quantum metric space for $ k >
2.$ \ethm {\it Proof:} In view of theorem \ref{P:thm:7.4} and the
previous lemma we only have to show that $ (\cla_\nu, L_k)$ has
finite radius. Let $\mu_1,\mu_2 \in S(\cla_\nu) , a \in \cla_\nu $
with $ L_k (a) \leq 1 $. By lemma \ref{M:lem:2.2} $ \| a \| \leq C
$, because $ L_2(a) \leq L_k(a)$. Hence $ | \mu_1 (a ) - \mu_2 ( a
) | \leq 2 C$, that is $ diam (\cla_\nu, L_k) \leq 2 C.$ \qed
\bppsn \label{metric:ppsn:2.5} Let
\begin{displaymath}
0\longrightarrow A_0 \stackrel i\longrightarrow  A_1 \stackrel\pi
\longrightarrow A_2 \longrightarrow 0
\end{displaymath}
be a short exact sequence of $C^*$-algebras, with $A_1,A_2$ unital
and a positive linear splitting $ \sigma : A_2 \raro A_1 $.  Let $
\phi : A_1^\prime  \raro A_0^\prime \oplus   A_2^\prime, \psi :
A_0^\prime \oplus   A_2^\prime  \raro   A_1^\prime$ be the bounded
linear maps given by $$ \phi ( \mu ) = ( \mu_1,\mu_2), \mu_1=  \mu
|_{i (A_0)}       ,\mu_2   =  \mu \circ \sigma$$   $$ \psi  (
\mu_1,\mu_2)= \mu, \mu (a ) = \mu_2 ( \pi ( a)) + \mu_1 ( a -
\sigma \circ \pi ( a )) $$ Then $ \mu_1, \mu_ 2$ are inverse to
each other. \eppsn {\it Proof:} Let $  \phi ( \mu ) = (
\mu_1,\mu_2),  \psi  (  \mu_1,\mu_2) = \mu^\prime $. Then \bean
\mu^\prime ( a ) & = & \mu_2 ( \pi ( a)) + \mu_1 ( a - \sigma
\circ \pi ( a )) \cr & = & \mu ( \sigma \circ \pi ( a)) + \mu ( a
- \sigma \circ \pi ( a )) \cr & = & \mu ( a ). \eean Therefore $
\psi \circ \phi = {Id}_ { A_1^\prime }$. Similarly one can show
that the other composition is also identity. \qed\\
Let $\cla,Lip(\cla),L$ be as above. Suppose we have    a short
exact sequence of $C^*$-algebras
\begin{displaymath}
0\longrightarrow \clk \otimes \cla  \stackrel i\longrightarrow
\widetilde {\cla_1} \stackrel\pi \longrightarrow \widetilde
{\cla_2} \longrightarrow 0
\end{displaymath}
  with $\widetilde {\cla_1},\widetilde {\cla_2} $ unital and a positive unital linear
   splitting $ \sigma : \widetilde {\cla_2} \raro \widetilde {\cla_1}$. Let $(\cla_2, L_2) $
    be a compact quantum metric space with $\cla_2$ a dense subspace of selfadjoint elements
     of $\widetilde {\cla_2}$.  Define $\cla_1 = i (\widetilde {\cla_\nu}  ) \oplus \sigma
     ( \cla_2)$. Then we have
\bthm In the above set up $ L_1: \cla_1 \raro \IR_{+}$, given by
$$ L_1 ( a ) = L_2  ( \pi (( a ) ) + L_k  (  a  - \sigma \circ \pi
( a )) $$ is a Lip norm for $k > 2$. \ethm
{\it Proof:} We break the proof in several steps.\\
{\bf Step (i)} $L_1(a)=0$ iff $ a \in \IR {Id}_{A_1} $: If part is
obvious for the only if part note $L_1(a)=0$ gives $ \pi (a ) =
\lambda {Id}_{A_2}
 $
 for some $ \lambda \in \IR$ and $ L_0( a -\lambda {Id}_{A_1}) =0$. Hence $ a = \lambda
 {Id}_{A_1}$.\\
{\bf Step (ii)} $ (A_1,L_1) $ has finite radius: Let $\mu,\lambda
\in S(\cla_1)$ and
$(\mu_1,\mu_2)=\phi(\mu),(\lambda_1,\lambda_2)=\phi(\lambda)$,
where $\phi$ is as in proposition \ref {metric:ppsn:2.5}. Then
from the norm estimate of $\phi$ obtained in proposition
\ref{metric:ppsn:2.5} we get $\|\mu_i\|,\|\lambda_i\| \le (1+
\|\sigma\|),$ for $i=1,2$ and positivity of $\sigma$ implies
$\|\mu_2\|=\|\lambda_2\|=1$. Let $x \in \cla_1$ with $L(x) \le 1$,
then \bean |\mu(x)-\lambda(x)| & = & |\mu_2(\pi(x))+\mu_1(x -
\sigma \circ \pi (x))- \lambda_2(\pi(x))-\lambda_1(x - \sigma
\circ \pi (x))| \cr & \le &
|\mu_2(\pi(x))-\lambda_2(\pi(x))|+|\mu_1(x - \sigma \circ \pi
(x))-\lambda_1(x - \sigma \circ \pi (x))| \cr & \le &
diam(\cla_2,L_2)+2 ( 1 + \| \sigma \|)C \eean where $C$ is the
constant obtained in lemma \ref{M:lem:2.2}. This proves
$(\cla_1,L_1)$
 has finite radius.\\
{\bf Step (iii)} In view of theorem \ref{P:thm:7.4} it suffices to
show that $\clb_1=\{ a \in \cla_1 : \; \| a \| \le 1, L_1(a) \le
1\}$ is totally bounded. Since $(\cla_\nu,L_k)$ and $(\cla_2,L_2)$
are compact quantum metric spaces it follows that if we have a
sequence
 $a_n \in \clb_1$, then there exists a subsequence $a_{n_k}$ such that both $\pi(a_{n_k})$
 and $a_{n_k} - \sigma \circ \pi(a_{n_k})$ converges in norm. Hence $a_{n_k}$ converges in
 norm implying the totally boundedness. \qed\\
\newsection{Examples}
\bxmpl  \rm Let $\Omega$ be a strongly pseudoconvex domain in
$\IC^n$. Let $H^2(\del \Omega)$ be the closure in $L^2(\del
\Omega)$ of boundary values of holomorphic functions that can be
continuously extended to $\bar{\Omega}$. For $f \in C(\del
\Omega)$ let $T_f$ be the associated Toeplitz operator, that is
the compression of the multiplication operator $M_f$ on $L^2(\del
\Omega)$ on $H^2(\del \Omega)$. Let $\frak{T}(\del \Omega)$ be the
associated Toeplitz extension, that is the $C^*$-algebra generated
by the operators $T_f$ along with the compacts. Then we have a
short exact sequence of $C^*$-algebras
\begin{displaymath}
0\longrightarrow \clk(H^2(\del \Omega)) \stackrel i\longrightarrow
\frak{T}(\del\Omega) \stackrel\pi \longrightarrow C(\del \Omega)
\longrightarrow 0
\end{displaymath}
Since this sequence admits a positive unital splitting by the
previous theorem we get CQMS structure on $\frak{T}(\del \Omega)$.
\exmpl \bxmpl  \rm The $C^*$-algebra of continuous functions on
the quantum $SU(2)$, to be denoted by $C(SU_q(2))$, is the
universal $C^*$-algebra generated by two elements $\alpha$ and
$\beta$ satisfying the following relations: \bean
\alpha^*\alpha+\beta^*\beta=I,&&\alpha\alpha^*
+q^2\beta\beta^*=I,\\
\alpha\beta-q\beta\alpha=0,&&\alpha\beta^*-q\beta^*\alpha=0,\\
\beta^*\beta&=&\beta\beta^*. \eean
 The $C^*$-algebra
$C(SU_q(2))$ can be described more concretely as follows. Let
$\{e_i\}_{i\geq0}$ and $\{e_i\}_{i\in\IZ}$ be the canonical
orthonormal bases for $L_2(\IN_0)$ and $L_2(\IZ)$ respectively. We
denote by the same symbol $N$ the operator $e_k\mapsto ke_k$,
$k\geq0$, on $L_2(\IN_0)$ and $e_k\mapsto ke_k$, $k\in\IZ$, on
$L_2(\IZ)$. Similarly, denote by the same symbol $\ell$ the
operator $e_k\mapsto e_{k-1}$, $k\geq1$, $e_0\mapsto0$ on
$L_2(\IN_0)$ and the operator $e_k\mapsto e_{k-1}$, $k\in\IZ$ on
$L_2(\IZ)$. Now take $\clh$ to be the Hilbert space
$L_2(\IN_0)\ots L_2(\IZ)$, and define $\pi$ to be the following
representation of $C(SU_q(2))$ on $\clh$:
\[
\pi(\alpha)=\ell\sqrt{I-q^{2N}}\ots I,\;\;\;\pi(\beta)=q^N\ots
\ell.
\]
Then $\pi$ is a faithful representation of $C(SU_q(2))$, so that
one can identify $C(SU_q(2))$ with the $C^*$-subalgebra of
$\clb(\clh)$ generated by $\pi(\alpha)$ and $\pi(\beta)$. Image of
$\pi$ contains $\clk \otimes C( \IT )$ as an ideal with $C( \IT )$
as the quotient algebra, that is we have a useful short exact
sequence

\be \label{100} 0\longrightarrow \clk\otimes C( \IT ) \stackrel
i\longrightarrow  \cla \stackrel\sigma \longrightarrow C( \IT )
\longrightarrow 0. \ee The homomorphism $\sigma$ is explicitly
given by $\sigma(\alpha)=\ell, \sigma(\beta)=0$.  It is easy to
see that the above short exact sequence admits a positive
splitting taking $z^n \in C( \IT )$ to ${\ell}^n\otimes I$, for
all $n \ge 0$. Hence we get a compact quantum metric space
structure on $C(SU_q(2))$. \exmpl
 \bxmpl  \rm
Quantum sphere was introduced by Podles in \cite {PO}. This is the
universal $C^*$-algebra denoted by $C(S^2_{qc})$, generated by two
elements $A$ and $B$ subject to the following relations: \bean
A^* = A ,&&  B^*B = A - A^2 + c I ,\\
BA = q^2 AB,&&BB^* = q^2 A - q ^ 4 + c I. \eean Here the
deformation parameters $ q,c$ satisfy $ | q | < 1, c>0 $. For
later purpose we also note down two irreducible representations
such that the representation given by the direct sum of these two
is faithful. Let $\clh_{+}  = l^2 ( \IN_0  ), \clh_{-} = \clh_{+}.
$ Define $ \pi_{\pm} (A), \pi_{\pm} (B): \clh_\pm \raro \clh_\pm$
by \bean \pi_{\pm}(A) ( e_n)= \lambda_\pm q^{2n} e_n \;\;\quad \;
&{\rm   where   }& \quad \; \lambda_\pm = \frac {1}{ 2}  \pm {(c +
\frac {1}{4} )}^{1/2} \cr \pi_{\pm}(B) ( e_n)= {c_\pm (n)}^{1/2}
e_{n-1} &{\rm  where}&
 c_\pm (n)= \lambda_\pm q^{2n}- {(\lambda_\pm q^{2n})}^2 + c , {\mbox{and}} e_{-1}=0 .
\eean
 Since $\pi =\pi_+\oplus \pi_-$ is a faithful representation  we
 have (\cite{SH}),\\
(i)$ C (S^2_{qc}) \cong C^* (\frak{T})\oplus_\sigma C^*(\frak{T}) := \{(x,y) : \; x,y \in C^*(\frak{T}), \sigma(x)=\sigma(y) \} $ where $C^*(\frak{T})$ is the Toeplitz algebra and $\sigma : C^*(\frak{T}) \raro C( \IT )$ is the symbol homomorphism.\\
(ii) We have a short exact sequence \bea \label{S:ses:2}
0\longrightarrow \clk \stackrel i\longrightarrow  C(S^2_{qc})
\stackrel\alpha \longrightarrow C^*(\frak {T}) \longrightarrow 0
\eea
 As in the earlier case this
short exact sequence  is also split exact. Here a positive
splitting is given by $\ell \in C^*(\frak {T}) \mapsto
(\ell,\ell)$. Now to apply the basic theorem note that by the
earlier example on Toeplitz extensions we already have a Lip norm
on a dense subspace of $C^*(\frak {T})$. \exmpl

\end{document}